\documentclass[10pt,twoside]{amsart}
\usepackage{lmodern}
\usepackage[T1]{fontenc}
\usepackage{amsmath,amsfonts,amsthm,amssymb,amscd}
\usepackage{mathrsfs}
\usepackage{latexsym}
\usepackage[pdftex]{graphicx}
\usepackage{enumitem}
\usepackage{bbm}
\usepackage{accents}

\usepackage{geometry}
\newtheorem{corollary}{Corollary}[section]
\newtheorem{lemma}{Lemma}[section]
\newtheorem{theorem}{Theorem}[section]

\usepackage{hyperref}  %added
\hypersetup{colorlinks=true, linktoc=all,
    linkcolor=blue,
		citecolor=blue,
		urlcolor=cyan,
		bookmarksopen=true
		}
\usepackage{hyperref}  %added
\newcommand{\lnc}{\mathscr{L}}

\title{Sum of short exponential sums with prime numbers}
\author{Firuz Rakhmonov}
\address{A.Dzhuraev Institute of Mathematics,  National Academy of Sciences of Tajikistan}
\email{rakhmonov.firuz@gmail.com}
\date{}

\begin{document}

\begin{abstract}
For sufficiently large integers $K$, $x$, $y$, and $q$ satisfying $K \le y < x$, where
$f(u) = \alpha u^n + \alpha_{n-1}u^{n-1} + \ldots + \alpha_1 u$
is a polynomial of degree $n$ with real coefficients,
$n$ is a fixed positive integer,
$\alpha$ is a real number such that
$\left|\alpha - \frac{a}{q}\right| \le \frac{1}{q^2}$,
$(a, q) = 1$, $q \ge 1$ and $\lnc  = \ln x$,
an estimate of the form
$$
\sum_{k=1}^K \left| \sum_{x - y < p \le x} e(kf(p)) \right|
\ll K y \left( \frac{1}{q} + \frac{1}{y} + \frac{q}{K y^n} + \frac{1}{K^{2^{n-1}}} \right)^{2^{-n-1}} {\lnc}^{\frac{n^2}{2^{n+1}}},
$$
is obtained, which represents a strengthening and generalization of the corresponding estimate of I.~M.~Vinogradov.


\emph{Bibliography:} 18 references.
\end{abstract}

\maketitle

\section{Introduction}
The method of estimating exponential sums with prime numbers developed by I.~M.~Vinogradov enabled him to solve a number of arithmetic problems involving primes.
One of these problems concerns the distribution of the fractional parts $\{\alpha p\}$.
In this direction, I.~M.~Vinogradov~\cite{VinigradovIM-Izb-Trud.eng,Vinogradov+Karatsuba-TrMIAN-1984-77.eng} obtained a considerably better estimate for an exponential sum than in the general case:

\emph{Let $K$ and $x$ be sufficiently large integers with $K \le x$, and let $\alpha$ be a real number of the form}
\begin{equation}\label{formula opr alpha}
\alpha = \frac{a}{q} + \frac{\theta}{q^2}, \qquad (a,q)=1, \quad q \ge 1, \qquad |\theta| \le 1,
\end{equation}
\emph{then, for $q \le x$, we have}
\[
V_1(K,x) = \sum_{k=1}^K \left| \sum_{p \le x} e(\alpha k p) \right|
\ll K x^{1+\varepsilon} \left( \sqrt{ \frac{1}{q} + \frac{q}{x} } + x^{-0.2} \right).
\]

It follows that the following statement holds:
\emph{For any $\sigma$ satisfying $0 < \sigma \le 1$, the number $F_\alpha(x,\sigma)$ of values of $\{\alpha p\}$, $p \le x$, satisfying the condition $\{\alpha p\} < \sigma$, is given by the formula}
\[
F_\alpha(x,\sigma) = \sigma \pi(x) + R_\sigma(x),
\qquad
R_\sigma(x) \ll x^{1+\varepsilon} \left( \sqrt{ \frac{1}{q} + \frac{q}{x} } + x^{-0.2} \right).
\]
In particular, if $\alpha$ is an irrational number with bounded partial quotients, one can choose $q$ to be of order $\sqrt{x}$.
In this case, in the problem of the distribution of the fractional parts $\{\alpha p\}$, where the argument runs over the prime numbers in a short interval, that is,
\emph{for $F_\alpha(x, y, \sigma)$ --- the number of terms of the sequence $\{\alpha p\}$ such that $x - y < p \le x$ and $\{\alpha p\} < \sigma$,
bearing in mind that $F_\alpha(x, y, \sigma) = F_\alpha(x, \sigma) - F_\alpha(x - y, \sigma)$, the following asymptotic formula holds}
\begin{equation*}
F_\alpha(x, y, \sigma) = \sigma \bigl( \pi(x) - \pi(x - y) \bigr) + O\bigl( x^{\frac{4}{5} + \varepsilon} \bigr),
\end{equation*}
which is nontrivial for $y \gg x^{\frac{4}{5} + \varepsilon}$.

Using the method of estimating exponential sums with prime numbers developed by I.~M.~Vinogradov~\cite{VinigradovIM-Izb-Trud.eng,Vinogradov+Karatsuba-TrMIAN-1984-77.eng},
in combination with the methods of~\cite{RakhmonovZKh-1993-IzvRAN.eng,RakhmonovFiruz-VestMGU-2011-garv},
it was proved in~\cite{Rakhmonovi-DoklMath2014_124.eng,RakmonovZKh+RakmonovFZ+IsmatovSN-DANRT2013.eng} that:
\emph{Let $K$, $x$, $y$, and $q$ be sufficiently large integers with $K \le y$, let $A$ be an absolute constant, $\lnc = \ln x$, and let $\alpha$ be a real number of the form \textup{(\ref{formula opr alpha})}.
Then, for $y \gg x^{\frac{2}{3}} (\lnc)^{4A + 16}$ and $\lnc^{4A + 20} \le q \le K y^2 x^{-1} (\lnc)^{-4A - 20}$, the following estimate holds:}
\[
V_1(K, x, y) = \sum_{k=1}^K \left| \sum_{x - y < p \le x} e(\alpha k p) \right|
\ll \frac{K y}{\lnc^A}.
\]

The main result of this paper is Theorem~\ref{Teorema sumKorTrigSummSPr},
which provides an estimate for the sums of moduli of short exponential sums with prime numbers.

\begin{theorem} \label{Teorema sumKorTrigSummSPr}
Let $K$, $x$, $y$, and $q$ be sufficiently large integers satisfying $K \le y < x$, and let $n$ be a fixed natural number.
Suppose that
$$
f(u) = \alpha u^n + \alpha_{n-1}u^{n-1} + \ldots + \alpha_1 u
$$
is a polynomial of degree $n$ with real coefficients, where
\begin{equation*}
\alpha = \frac{a}{q} + \frac{\theta}{q^2}, \qquad (a, q) = 1, \quad q \ge 1, \qquad |\theta| \le 1.
\end{equation*}
Then the following estimate holds:
$$
V_n(K,x,y) = \sum_{k=1}^K \left| \sum_{x - y < p \le x} e(kf(p)) \right|
\ll Ky \left( \frac{1}{q} + \frac{1}{y} + \frac{q}{K y^n} + \frac{1}{K^{2^{n-1}}} \right)^{2^{-n-1}} {\lnc}^{\frac{n^2}{2^{n+1}}},
$$
where $\lnc = \ln x$.
\end{theorem}

The proof of Theorem~\ref{Teorema sumKorTrigSummSPr} is based on Weyl's method combined with ideas developed in
\cite{RakhmonovZKh+FZ-ChebSbornik-2024-25-2.eng,RakhmonovFZ-2024-ChebSb-25-4(95).eng,RakhmonovFZ-2024-DNANT-67-5-6.eng,RZKh+RFZ-DNAT-2023-66-9-10.eng,RZKh+RFZ-DNAT-2024-67-3-4.eng}.
In these papers, when deriving asymptotic formulas for generalized Waring and Estermann problems with almost proportional terms,
Hua Loo-Keng's theorem (\cite{Vaughan-1985-Mir-Metod Hardy-Littlvuda.eng}, Lemma~2.5) on the mean value of the exponential sum
\[
T(\alpha, x) = \sum_{m \le x} e(\alpha m^n)
\]
was extended to short Weyl exponential sums of the form
\[
T(\alpha; x, y) = \sum_{x - y < m \le x} e(\alpha m^n),
\]
and a nontrivial estimate for $T(\alpha; x, y)$ on minor arcs was obtained (see also
\cite{RakhmonovFZ-ChebSbornik-2011-12-1.eng,RakhmonovFZ-VestnikMGU-2011-3.eng,RakhmonovZKh+RFZ-TrMIAN-296-2016.eng,RakhmonovZKh+FZ-ChebSbornik-2019-20-4.eng}).

We now determine the range of parameters $K$, $q$, and $y$ for which the estimate obtained in
Theorem~\ref{Teorema sumKorTrigSummSPr} is nontrivial.

\begin{corollary}\label{Sled1 teorema sumKorTrigSummSPr}
Assume that the conditions of Theorem~\ref{Teorema sumKorTrigSummSPr} are satisfied and that
\begin{equation*}
K \gg (\lnc)^{4A + \frac{n^2}{2^{\,n-1}}}, \qquad
(\lnc)^{2^{\,n+1}A + n^2} \ll q \ll K y^2 (\lnc)^{-2^{\,n+1}A - n^2}, \qquad
y \gg (\lnc)^{2^{\,n+1}A + n^2},
\end{equation*}
where $A \ge 2$ is an absolute fixed constant.
Then the following estimate holds:
\[
V_n(K, x, y) \ll K y \lnc^{-A}.
\]
\end{corollary}

\section{Auxiliary statements}

{\lemma \label{Lemma razn oper proizv mnog} Let $\Delta_j$ denote the $j$-th application of the difference operator, so that for any polynomial $f(u)=a_su^s+\ldots+a_1u$ we have
\begin{align}
&\Delta_1(f(u);h)=f(u+h)-f(u),
&\Delta_{j+1}(f(u);h_1,\ldots,h_{j+1})=\Delta_1(\Delta_s(f(u);h_1,\ldots,t_j);h_{j+1}). \label{form-lemm-razn-oper}
\end{align}
Then for $j=1,\ldots,s-1$ the following relation holds:
\begin{equation*}
\begin{split}
\Delta_j(f(u);h_1,\ldots,h_j)&=h_1h_2\ldots h_j\sum_{i_0=0}^{s-j}a_{s-i_0}\sum_{i_1=0}^{s-j-i_0}C_{s-i_0}^{i_1+1} \sum_{i_2=0}^{s-j-i_0-i_1} C_{s-1-i_0-i_1}^{i_2+1}\ldots\times\\ &\times\sum_{i_j=0}^{s-j-i_0-i_1-\ldots-i_{j-1}}C_{s-j+1-i_0-i_1-\ldots-i_{j-1}}^{i_j+1}h_1^{i_1}h_2^{i_2}\ldots h_j^{i_j} u^{s-j-i_0-i_1-\ldots-i_j}.
\end{split}
\end{equation*}}

{\sc Proof.} For convenience and symmetry, we use the identity
\begin{equation}\label{formula (u+h)r-ur=}
(u+h)^r-u^r
=h\sum_{i=0}^{r-1}C_r^{i+1}h^iu^{r-1-i}.
\end{equation}
Let us find $\Delta_1(f(u);h_1)$. Applying formula (\ref{form-lemm-razn-oper}) and then, for $r=s-i_0$, formula (\ref{formula (u+h)r-ur=}), we obtain
\begin{align*}
\Delta_1(f(u);h_1)=&\sum_{i_0=0}^{s-1}a_{s-i_0}((u+h_1)^{s-i_0}-u^{s-i_0})=
h_1\sum_{i_0=0}^{s-1}a_{s-i_0}\sum_{i_1=0}^{s-1-i_0}C_{s-i_0}^{i_1+1}h_1^{i_1}u^{s-1-i_0-i_1}.
\end{align*}
Using formula (\ref{form-lemm-razn-oper}) with $r=s-1-i_0-i_1$, we find $\Delta_2(f(u);h_1,h_2)$:
\begin{align*}
\Delta_2((f(u);h_1,h_2)&=\Delta_1\left(\Delta_1((f(u);h_1);h_2\right)=\Delta_1 \left(h_1\sum_{i_0=0}^{s-1}a_{s-i_0}\sum_{i_1=0}^{s-1-i_0}C_{s-i_0}^{i_1+1}h_1^{i_1}u^{s-1-i_0-i_1};h_2\right)=\\
&=h_1\sum_{i_0=0}^{s-1}a_{s-i_0}\sum_{i_1=0}^{s-1-i_0}C_{s-i_0}^{i_1+1}h_1^{i_1}\left((u+h_2)^{s-1-i_0-i_1}-u^{s-1-i_0-i_1}\right)=
\\ &=
 h_1h_2\sum_{i_0=0}^{s-1}a_{s-i_0}
\sum_{i_1=0}^{s-1-i_0}C_{s-i_0}^{i_1+1}\sum_{i_2=0}^{s-2-i_0-i_1}C_{s-1-i_0-i_1}^{i_2+1} h_1^{i_1}h_2^{i_2} u^{s-2-i_0-i_1-i_2}.
\end{align*}
From the upper limit of the sum over $i_2$ it follows that $i_0+i_1+i_2\le s-2$. Consequently, $i_1\le s-2-i_0$ and $i_0\le s-2$, so the right-hand side of the last equality takes the form
\begin{align*}
&\Delta_2((f(u);h_1,h_2)=h_1h_2\sum_{i_0=0}^{s-2}a_{s-i_0}
\sum_{i_1=0}^{s-2-i_0}C_{s-i_0}^{i_1+1}\sum_{i_2=0}^{s-2-i_0-i_1}C_{s-1-i_0-i_1}^{i_2+1} h_1^{i_1}h_2^{i_2} u^{s-2-i_0-i_1-i_2}.
\end{align*}
Assuming that the lemma holds for $j\le s-2$, let us prove it for $j+1$. Using relation (\ref{form-lemm-razn-oper}) and then, for $r=s-j-i_0-i_1-\ldots-i_j$, formula (\ref{formula (u+h)r-ur=}), we obtain
\begin{align*}
\Delta_{j+1}(f(u)&;h_1,\ldots,h_{j+1})=\Delta_1\left(\Delta_j(f(u);h_1,\ldots,h_j),h_{j+1}\right)=\\
=\Delta_1&\left(h_1h_2\ldots h_j\sum_{i_0=0}^{s-j}a_{s-i_0}\sum_{i_1=0}^{s-j-i_0}C_{s-i_0}^{i_1+1} \sum_{i_2=0}^{s-j-i_0-i_1} C_{s-1-i_0-i_1}^{i_2+1}\ldots\times\right.\\ &\left.\times\sum_{i_j=0}^{s-j-i_0-i_1-\ldots-i_{j-1}}C_{s-j+1-i_0-i_1-\ldots-i_{j-1}}^{i_j+1}h_1^{i_1}h_2^{i_2}\ldots h_j^{i_j} u^{s-j-i_0-i_1-\ldots-i_j};\ h_{j+1}\right)=\\
=h_1&h_2\ldots h_j\sum_{i_0=0}^{s-j}a_{s-i_0}\sum_{i_1=0}^{s-j-i_0}C_{s-i_0}^{i_1+1} \sum_{i_2=0}^{s-j-i_0-i_1} C_{s-1-i_0-i_1}^{i_2+1}\ldots\sum_{i_j=0}^{s-j-i_0-i_1-\ldots-i_{j-1}}\times\\
&\times C_{s-j+1-i_0-i_1-\ldots-i_{j-1}}^{i_j+1}h_1^{i_1}h_2^{i_2}\ldots h_j^{i_j} \left((u+h_{j+1})^{s-j-i_0-i_1-\ldots-i_j}-u^{s-j-i_0-i_1-\ldots-i_j}\right)=
\end{align*}
\begin{align*}
=h_1&\ldots h_{j+1}\sum_{i_0=0}^{s-j}a_{s-i_0}\sum_{i_1=0}^{s-j-i_0}C_{s-i_0}^{i_1+1} \sum_{i_2=0}^{s-j-i_0-i_1} C_{s-1-i_0-i_1}^{i_2+1}\ldots\sum_{i_j=0}^{s-j-i_0-i_1-\ldots-i_{j-1}}C_{s-j+1-i_0-i_1-\ldots-i_{j-1}}^{i_j+1}\times\\
&\times \sum_{i_{j+1}=0}^{s-j-1-i_0-i_1-\ldots-i_j}C_{s-j-1-i_0-i_1-\ldots-i_j}^{i_{j+1}+1}h_1^{i_1}h_2^{i_2}\ldots h_{j+1}^{i_{j+1}} u^{s-j-1-i_0-i_1-\ldots-i_{j+1}}.
\end{align*}
We proceed as in the case of $\Delta_2(f(u);h_1,h_2)$. From the upper limit of the sum over $i_{j+1}$ it follows that
$$
i_0+i_1+\ldots+i_j+i_{j+1}\le s-j-1,
$$
which in turn implies that
$$
i_j\le s-j-1-i_0-i_1-\ldots-i_{j-1},\qquad
i_1\le s-j-1-i_0, \qquad i_0\le s-j-1.
$$
Therefore, the right-hand side of the last equality takes the form
\begin{align*}
\Delta_{j+1}(f(u);h_1,\ldots,h_{j+1})= h_1\ldots h_{j+1}\sum_{i_0=0}^{s-j-1}a_{s-i_0}\sum_{i_1=0}^{s-j-1-i_0}C_{s-i_0}^{i_1+1} \sum_{i_2=0}^{s-j-i_0-i_1} C_{s-1-i_0-i_1}^{i_2+1}\ldots\times\\
\times \sum_{i_{j+1}=0}^{s-j-1-i_0-i_1-\ldots-i_j}C_{s-j-1-i_0-i_1-\ldots-i_j}^{i_{j+1}+1}h_1^{i_1}h_2^{i_2}\ldots h_{j+1}^{i_{j+1}} u^{s-j-1-i_0-i_1-\ldots-i_{j+1}}.
\end{align*}
Lemma \ref{Lemma razn oper proizv mnog} is proved.

{\lemma\label{Lemma virazh st T(f(u);x,y) konech razn}  {\rm (\cite{RakhmonovFZ-2024-ChebSb-25-4(95).eng}}.
Let $f(u)$ be a polynomial of degree $s$, and let $x$ and $y$ be positive integers such that $y<x$. Then for $j=1,\ldots,s-1$ the following inequality holds:
\begin{align*}
&\left|\sum_{x-y<u\le x}e(f(u))\right|^{2^{j}}\leq(2y)^{2^j-j-1}\sum_{|h_1|<y}\ldots\sum_{|h_j|<y}\left|\sum_{u\in I_j(x,y;h_1,\ldots,h_j)}e(\Delta_j(f(u); h_1,\ldots,h_j))\right|,
\end{align*}
where the intervals $I_j(x,y;h_1,\ldots,h_j)$ are defined by the relations:
\begin{align*}
&I_1(x,y;h_1)=(x-y,x]\cap (x-y-h_1,x-h_1],\\
&I_j(x,y;h_1,\ldots, h_{j})=I_{j-1}(x,y;h_1,\ldots, h_{j-1})\cap I_{j-1}(x-h_j,y;h_1,\ldots, h_{j-1}),
\end{align*}
that is, $I_{j-1}(x-h_j,y;h_1,\ldots, h_{j-1})$ is obtained from $I_{j-1}(x,y;h_1,\ldots, h_{j-1})$ by shifting by $-h_j$ of all intervals, which forms it as an intersection.}

{\lemma \label{Lemma ob otsenke summ st obsh funk delit} {\rm \cite{Mardjanishvili-garv}.} For $x\ge 1$, $r\ge2$, and $k\ge1$, we have
\begin{align*}
\sum_{h\le x}\tau_r^k(h)\ll\frac{xr^k}{(r!)^{\frac{r^k-1}{r-1}}}\left(\ln x+r^k-1\right)^{r^k-1}.
\end{align*}}

{\lemma \label{Lemma ob otsenke summ obr RDBTsCh} {\rm \cite{Karatsuba-OATCh.eng}.} Let $\alpha$ be a real number such that
$$
\left|\alpha-\frac{a}{q}\right|\le \frac{1}{q^2}, \qquad (a,q)=1,
$$
where $x\ge1$, $y>0$, and $\beta$ is arbitrary. Then the following estimate holds:
$$
\sum_{h\le x}\min\left(y,\ \frac1{\|\alpha n+\beta\|}\right)\le6\left(\frac{x}q+1\right)(y+q\ln q).
$$}

\section{Proof of Theorem \ref{Teorema sumKorTrigSummSPr}.}
Without loss of generality, we may assume that all coefficients of the polynomial
$$
f(u) = \alpha u^n + \alpha_{n-1}u^{n-1} + \ldots + \alpha_1u
$$
except for the leading one are equal to zero; that is, $f(u) = \alpha u^n$.
In this case, the sum $V_n(K,x,y)$ takes the following form:
$$
V_n(K,x,y) = \sum_{k=1}^K \left| \sum_{x - y < p \le x} e(\alpha k p^n) \right|.
$$

Applying the Cauchy inequality, taking the sum over $k$ as the inner one, and dividing the sum over $p_1$ and $p_2$ into three parts corresponding to the conditions $p_1<p_2$, $p_1=p_2$, and $p_1>p_2$, and estimating the sum under the condition $p_1=p_2$ by a quantity $\ll K^2y$, and also noting that the absolute values of the sums under the conditions $p_1<p_2$ and $p_1>p_2$ are equal, we obtain
\begin{align*}
V_n^2(K,x,y)&\le K\sum_{k=1}^K\sum_{x-y<p_1,p_2\le x}e(\alpha k(p_1^n-p_2^n))\le\\
&\le K^2y+2K\sum_{x-y<p_1<x}\sum_{x-y<p_2<p_1}\left|\sum_{k=1}^Ke(k(p_1^n-p_2^n))\right|.
\end{align*}
In the right-hand side of the last sum, replacing the summation over the prime numbers $p_1$ and $p_2$ by summation over natural numbers $m_1$ and $m_2$, and noting that the inequalities $x-y<m_2<m_1$ and $0<m_1-m_2<m_1-x+y$ are equivalent, and then introducing the notation $m_2=m_1-m$, as well as taking into account that
\begin{equation}\label{formula opr mnogochelena f(m1)}
\begin{split}
  & m_1^n-m_2^n=m_1^n-(m_1-m)^n=mf(m_1),\\
  &f(m_1)=a_{n-1}m_1^{n-1}+\ldots+a_1m_1,\qquad a_{n-i}=(-1)^{i-1}C_n^im^{i-1},
\end{split}
\end{equation}
(that is, $f(m_1)$ is a polynomial of degree $n-1$ in $m_1$), we consequently find
\begin{align}
V_K^2&(x,y)\le K^2y+2K\sum_{x-y<m_1<x}\sum_{x-y<m_2<m_1}\left|\sum_{k=1}^Ke(\alpha k(m_1^n-m_2^n))\right|=\nonumber\\
&=K^2y+K\sum_{x-y<m_1<x}\sum_{0<m_1-m_2<m_1-x+y}\left|\sum_{k=1}^Ke(\alpha k(m_1^n-m_2^n))\right|=\nonumber\\
&=K^2y+K\sum_{x-y<m_1<x}\sum_{0<m<m_1-x+y}\left|\sum_{k=1}^Ke(\alpha kmf(m_1))\right|=\nonumber\\
&=K^2y+KW_1,\quad  W_1=\sum_{0<m<y}\sum_{x-y+m<m_1<x}\left|\sum_{k=1}^Ke(\alpha kmf(m_1))\right|. \label{formula VK2 cherez W_1}
\end{align}

Raising $W_1$ to the square and applying the Cauchy inequality twice, we have
\begin{align*}
W_1^2&\le y^2\sum_{0<m<y}\sum_{x-y+m<u<x}\left|\sum_{k=1}^Ke(\alpha kmf(u))\right|^2=\\
&=y^2\sum_{0<m<y}\sum_{x-y+m<u<x}\sum_{k_1=1}^K\sum_{k_2=1}^Ke(\alpha(k_1-k_2)mf(u)),
\end{align*}
Dividing the sum over $k_1$ and $k_2$ into three parts corresponding to the conditions $k_1<k_2$, $k_1=k_2$, and $k_1>k_2$, and estimating the sum under the condition $k_1=k_2$ by a quantity $\ll y^4K$, and also noting that the absolute values of the sums under the conditions $k_1<k_2$ and $k_1>k_2$ are equal, we obtain
\begin{align}
&W_1^2\ll y^4K+y^2W_2, &W_2=\sum_{0<m<y}\sum_{x-y+m<u<x}\sum_{2\le k_1\le K}\sum_{1\le k_2<k_1}e(\alpha(k_1-k_2)mf(u)).\label{formula W1 cherez W2}
\end{align}
Using the fact that the inequalities $1\le k_2<k_1$ and $1\le k_1-k_2\le k_1-1$ are equivalent, introducing the notation $k=k_1-k_2$ and taking the sum over $u$ as the inner one, and then proceeding to estimates, we consecutively obtain
\begin{align*}
W_2&=\sum_{0<m<y}\sum_{x-y+m<u<x}\sum_{2\le k_1\le K}\sum_{1\le k_1-k_2\le k_1-1}e(\alpha(k_1-k_2)mf(u))=\\
&=\sum_{0<m<y}\sum_{x-y+m<u<x}\sum_{2\le k_1\le K}\sum_{1\le k\le k_1-1}e(\alpha kmf(u))=\\
&=\sum_{0<m<y}\sum_{2\le k_1\le K}\sum_{1\le k\le k_1-1}\sum_{x-y+m<u<x}e(\alpha kmf(u))\le\\
&\le K\sum_{0<m<y}\sum_{k\le K}\left|\sum_{x-y+m<u<x}e(\alpha kmf(u))\right|.
\end{align*}
Raising both sides of the last inequality to the power $2^{n-2}$ and applying the Cauchy inequality twice, and then applying Lemma \ref{Lemma virazh st T(f(u);x,y) konech razn} to the inner sum over $u$ for $s=n-1$ and $j=n-2$, we have
\begin{align}
W_2^{2^{n-2}}\le K^{2^{n-1}-1}y^{2^{n-2}-1}\sum_{m<y}\sum_{k\le K}\left|\sum_{x-y+m<u<x}e(\alpha kmf(u))\right|^{2^{n-2}}\ll  K^{2^{n-1}-1}y^{2^{n-1}-n}\cdot W_3\nonumber\\
W_3=\sum_{k\le K}\sum_{m<y}\sum_{|h_1|<y-m}\ldots\sum_{|h_{n-2}|<y-m}\left|\sum_{u\in I_{n-2}(x,y-m;h_1,\ldots,h_{n-2})}\hspace{-30pt}e(\Delta_{n-2}(\alpha kmf(u); h_1,\ldots,h_{n-2}))\right|.\label{formula W2<=}
\end{align}
Using the linearity property of the difference operator and then Lemma \ref{Lemma razn oper proizv mnog} for $s=n-1$ and $j=n-2$, we have
\begin{align*}
\Delta_{n-2}(\alpha&kmf(u); h_1,\ldots,h_{n-2}))=\alpha km\Delta_{n-2}(f(u); h_1,\ldots,h_{n-2}))=\\
=&\alpha kmh_1h_2\ldots h_{n-2}\sum_{i_0=0}^1a_{n-1-i_0}\sum_{i_1=0}^{1-i_0}C_{n-1-i_0}^{i_1+1}\sum_{i_2=0}^{1-i_0-i_1} C_{n-2-i_0-i_1}^{i_2+1}\ldots\times\\
&\times\sum_{i_{n-2}=0}^{1-i_0-i_1-\ldots-i_{n-3}} C_{2-i_0-i_1-\ldots-i_{n-3}}^{i_{n-2}+1}h_1^{i_1}h_2^{i_2}\ldots h_{n-2}^{i_{n-2}} u^{1-i_0-i_1-\ldots-i_{n-2}}.
\end{align*}
The multiple sum over the summation variables $i_0,\ i_1,\ i_2,\ \ldots,\ i_{n-2}$ on the right-hand side of the last formula is a linear polynomial in $u$ and $h_1,\ldots,h_{n-2}$. In this sum, only one term is a linear monomial in $u$, which appears when all summation variables $i_0,\ i_1,\ i_2,\ \ldots,\ i_{n-2}$ are equal to zero; all other terms form a linear polynomial in $h_1,\ldots,h_{n-2}$ with integer coefficients depending on the parameter $n$, which we denote by $g(h_1,\ldots,h_{n-2})$. Extracting the linear monomial in $u$ and using the explicit value of the leading coefficient $a_{n-1}=n-1$ of the polynomial $f(u)$ from relation (\ref{formula opr mnogochelena f(m1)}), we have
\begin{align*}
\Delta_{n-2}&(\alpha kmf(u); h_1,\ldots,h_{n-2}))=\\
&=\alpha kmh_1h_2\ldots h_{n-2}\left(a_{n-1}\ C_{n-1}^1C_{n-2}^1\ldots C_2^1 u+g(h_1,\ldots,h_{n-2})\right)=\\
&=\alpha(n-1)!kmh_1h_2\ldots h_{n-2}u+\alpha kmh_1h_2\ldots h_{n-2}g(h_1,\ldots,h_{n-2})
\end{align*}

Substituting the right-hand side of the obtained equality into (\ref{formula W2<=}), we find
\begin{align*}
&W_3=\sum_{k\le K}\sum_{m<y}\sum_{|h_1|<y-m}\ldots\sum_{|h_{n-2}|<y-m}\left|\sum_{u\in I_{n-2}(x,y-m;h_1,\ldots,h_{n-2})}e\left(\alpha(n-1)!kmh_1h_2\ldots h_{n-2}u\right)\right|,
\end{align*}
In the last sum over $u$, the number of terms satisfying the relation $h_1\cdot\cdot\cdot h_{n-2}=0$ does not exceed $Ky\cdot(n-2)y(2y)^{n-3}$. Therefore,
\begin{align}\label{formula W2 cherez W4}
W_3&\ll 2^{n-2}W_4+(n-2)Ky^2(2y)^{n-3},\\
W_4=&\sum_{k\le K}\sum_{m<y}\sum_{1\le h_1<y-m}\ldots\sum_{1\le h_{n-2}<y-m}\left|\sum_{u\in I_{n-1}(x,y-m;h_1,\ldots,h_{n-1})}e\left(\alpha(n-1)!kmh_1h_2\ldots h_{n-2}u\right)\right|.\nonumber
\end{align}
The inner sum over $u$ is linear, and from the definition of the interval  $I_{n-1}(x,y-m;h_1,\ldots,h_{n-1})$ in Lemma \ref{Lemma virazh st T(f(u);x,y) konech razn}, it follows that
$$
\left|I_{n-1}(x,y-m;h_1,\ldots,h_{n-1})\right|\le y-m.
$$
Therefore,
\begin{align*}
W_4&\le\sum_{k\le K}\sum_{m<y}\sum_{1\le h_1<y-m}\ldots\sum_{1\le h_{n-2}<y-m}\min\left(y-m,\ \frac1{2\|\alpha(n-1)!kmh_1h_2\ldots h_{n-2}\|}\right)\le\\
&\le \sum_{h\le(n-1)!Ky^{n-1}}\tau_n(h)\min\left(y,\ \frac1{2\|\alpha h\|}\right).
\end{align*}
Applying the Cauchy inequality to the last sum, and then Lemmas 3 and 4, we obtain
\begin{align*}
W_4^2&\le y\sum_{h\le(n-1)!Ky^{n-1}}\tau_n^2(h) \sum_{h\le(n-1)!Ky^{n-1}}\min\left(y,\ \frac1{2\|\alpha h\|}\right)\ll\\
&\ll (n-1)!Ky^n\lnc^{n^2-1}\left(\frac{(n-1)!Ky^{n-1}}q+1\right)(y+q\ln q)\ll\\
&\ll K^2y^{2n}\left(\frac1q+\frac1{Ky^{n-1}}\right)\left(1+\frac qy\right)\lnc^{n^2}\ll\\
&\ll K^2y^{2n}\left(\frac1q+\frac1{Ky^{n-1}}+\frac1y+\frac{q}{Ky^n}\right)\lnc^{n^2}.
\end{align*}
Using formulas (\ref{formula W2 cherez W4}), (\ref{formula W1 cherez W2}), and (\ref{formula VK2 cherez W_1}), we express the sum $V_K(x,y)$ in terms of the sum $W_4$, and then substituting the obtained estimate for $W_4^2$, we get
\begin{align*}
V_K^{2^{n+1}}(x,y)&\ll K^{2^{n+1}}y^{2^{n+1}}\left(\frac1{y^{2^n}}+\frac1{K^{2^{n-1}}}+\frac1{y^2}+\frac{W_4^2}{K^2y^{2n}}\right)\ll\\
&\ll  K^{2^{n+1}}y^{2^{n+1}}\left(\frac1q+\frac1{Ky^{n-1}}+\frac1y+\frac{q}{Ky^n}+\frac1{K^{2^{n-1}}}\right)\lnc^{n^2}\ll\\
&\ll  K^{2^{n+1}}y^{2^{n+1}}\left(\frac1q+\frac1y+\frac{q}{Ky^n}+\frac1{K^{2^{n-1}}}\right)\lnc^{n^2}.
\end{align*}
Extracting the $2^{n+1}$-th root, we obtain the statement of Theorem \ref{Teorema sumKorTrigSummSPr}.


\begin{thebibliography}{12} % refs
\bibitem{VinigradovIM-Izb-Trud.eng} Vinogradov, I.~M., 1952, \emph{Izbrannye trudy. (Russian) [Selected works.]}, Izdat. Akad. Nauk SSSR, Moscow.
\bibitem{Vinogradov+Karatsuba-TrMIAN-1984-77.eng} Vinogradov, I.~M., \& Karatsuba,~A.~A., 1986, ``The method of trigonometric sums in number theory'', \emph{Proc. Steklov Inst. Math.}, vol.~168, pp.~3-30.
\bibitem{RakhmonovZKh-1993-IzvRAN.eng} Rakhmonov,~Z.~Kh., 1994,~''Theorem on the mean value of $\psi(x,\chi)$ and its applications'', \emph{Russian Academy of Sciences. Izvestiya Mathematics}, vol.~43, Is.~1, pp.~49~--~64.
\bibitem{RakhmonovFiruz-VestMGU-2011-garv} Rahmonov,~F.~Z., 2011, ''Estimate of quadratic trigonometric sums with prime numbers'', \emph{Vestnik Moskov. Univ. Ser. 1. Mat. Mekh.}, no.~3, pp.~56~--~60.
\bibitem{Rakhmonovi-DoklMath2014_124.eng} Rakhmonov,~Z.~Kh.,\& Rakhmonov,~F.~Z., 2014,~``Sum of short exponential sums over prime numbers'', \emph{Doklady Mathematics}, vol.~90, No~3, pp.~699~--~700.
\bibitem{RakmonovZKh+RakmonovFZ+IsmatovSN-DANRT2013.eng} Rakhmonov,~Z.~Kh., \& Rakhmonov,~F.~Z., Ismatov~S.~N.,~2013, ``Estimate of sums of short exponential sums over prime numbers'', \emph{Doklady Akademii nauk Respubliki Tadzhikistan}, vol.~56, no~12, pp.~937~--~945, (in Russian).
\bibitem{RakhmonovZKh+FZ-ChebSbornik-2024-25-2.eng}  Rakhmonov,~Z.~Kh.,\& Rakhmonov,~F.~Z., 2024,~``Asymptotic formula in the Waring's problem with almost proportional summands'', \emph{Chebyshevskii Sbornik}, vol.~25, Is.~2, pp.~138~--~168, (in Russian).
\bibitem{RakhmonovFZ-2024-ChebSb-25-4(95).eng} Rakhmonov,~F.~Z., 2024,~``Asymptotic formula in generalization of ternary Esterman problem with almost proportional summands'', \emph{Chebyshevskii Sbornik}, vol.~25, Is.~4, pp.~120~--~137, (in Russian).
\bibitem{RakhmonovFZ-2024-DNANT-67-5-6.eng} Rakhmonov,~F.~Z., 2024,~``Estimate of short G.Weyl exponential sums on minor arcs'', \emph{Doklady Natsional'noy akademii nauk Tadzhikistana}, vol.~67, Is~5-6,  pp.~238~--~242.
\bibitem{RZKh+RFZ-DNAT-2023-66-9-10.eng} Rakhmonov,~Z.~Kh.,\& Rakhmonov,~F.~Z., 2023,~``Waring's problem with almost proportional summands'',  \emph{Doklady Natsional'noy akademii nauk Tadzhikistana}, vol~66, Is~9-10, pp.~481~--~488.
\bibitem{RZKh+RFZ-DNAT-2024-67-3-4.eng} Rakhmonov,~Z.~Kh.,\& Rakhmonov,~F.~Z., 2024,~``Asymptotic formula in the Waring's problem with almost proportional summands'',  \emph{Doklady Natsional'noy akademii nauk Tadzhikistana}, vol~~67, Is~3-4,  pp.~125~--~136.
\bibitem{Vaughan-1985-Mir-Metod Hardy-Littlvuda.eng} Vaughan~R.~C., 1981. \emph{The Hardy-Littlewood method}, Cambridge Tracts in Mathematics, vol. 80, Cambridge University Press, Cambridge, 172~p.
\bibitem{RakhmonovFZ-ChebSbornik-2011-12-1.eng} Rakhmonov,~F.~Z., 2011,~``Estimation of trigonometric sums with prime numbers'', \emph{Chebyshevskii Sbornik}, vol.~12, Is.~1, pp.~158~--~171, (in Russian).
\bibitem{RakhmonovFZ-VestnikMGU-2011-3.eng} Rakhmonov,~F.~Z., 2011,~``Estimate of quadratic trigonometric sums with prime numbers'', \emph{Vestnik Moskovskogo Universiteta. Seriya 1. Matematika. Mekhanika},  Is.~3, pp.~56~--~60.
\bibitem{RakhmonovZKh+RFZ-TrMIAN-296-2016.eng} Rakhmonov,~Z.~Kh.,\& Rakhmonov,~F.~Z., 2017,~``Short Cubic Exponential Sums over Primes'',  \emph{Proceedings of the Steklov Institute of Mathematics}, vol.~296, pp.~211~--~233.
\bibitem{RakhmonovZKh+FZ-ChebSbornik-2019-20-4.eng} Rakhmonov,~Z.~Kh.,\& Rakhmonov,~F.~Z., 2019,~`` Trigonometric sums with the M\"obius function'',  \emph{Chebyshevskii Sbornik}, vol.~20, Is.~4, pp.~281~--~305, (in Russian).
\bibitem{Mardjanishvili-garv} Mardjhanashvili,~K.~K., 1939, ``An estimate for an arithmetic sum'', \emph{Doklady Akad. Nauk SSSR}, vol.~22, no~7, pp.~391~--~393.
\bibitem{Karatsuba-OATCh.eng} 	Karatsuba~A.~A., 1993, \emph{Basic analytic number theory}, Springer-Verlag, Berlin, xiv+222~pp.

\end{thebibliography}
\end{document}